\documentclass[12pt]{article}
\usepackage{amsfonts}
\usepackage{amsmath,latexsym,amssymb,amsfonts,amsbsy, amsthm}
\newtheorem{rmk}{Remark}[section]

\renewcommand{\theequation}{\thesection.\arabic{equation}}

\newtheorem*{theorem*}{Theorem}
\newtheorem*{definition*}{Definition}

\newtheorem{theorem}{Theorem}[section]
\newtheorem{lemma}[theorem]{Lemma}

\newtheorem{Theorem}{Theorem}[section]

\newtheorem{Proposition}[Theorem]{Proposition}

\begin{document}

\title{Remarks on the Prandtl Equation }
\author{{ DING Yutao}}
\date{Institute of Mathematics, AMSS, Academia Sinica,
Beijing; Hua Loo-Keng Key Laboratory of Mathematics, Chinese Academy
of Sciences.}

%\date{}
\maketitle

\begin{abstract}
In a recent result of G\'{e}rard-Varet and Dormy \cite{gd},  they established
ill-posedness for the Cauchy problem of the linearized Prandtl equation around non-monotic
special solution which is independent of $x$ and satisfies the heat equation. In \cite{gd2} and
\cite{gt}, some nonlinear ill-posedness were established with this counterexample.  Then it is
natural to consider the problem that does this linear ill-posedness happen whenever the  non-degenerate
critical points appear. In this paper, we concern the linearized Prandtl equation around general stationary
solutions with non-degenerate
critical points depending on $x$ which could be considered as the time-periodic solutions and show some ill-posdness.
\end{abstract}

\noindent {\sl Keywords:}  Prandtl Equations,

\vskip 0.2cm

%\noindent {\sl AMS Subject Classification (2000):} 35Q30, 76D03  \

\renewcommand{\theequation}{\thesection.\arabic{equation}}
\setcounter{equation}{0}
%%%%%%%%%%%%%%%%%%%%%%%%%%%%%%%%%%%%%%%%%%%%%%
%%%%%%%%%%%%%%%%%%%%%%%%%%%%%%%%%%%%%%%%%

\section{Introduction}

The behavior of the solution to the vanishing viscosity limit of Navier-Stokes equation
near a solid boundary is an outstanding open problem both in fluid mechanics and in mathematics.
To describe this problem, let us consider the two-dimensional Navier-Stokes equation on a half-space:
\begin{equation}
\begin{split}
 \partial_tu^{\nu}+u^\nu \partial_xu^\nu+v^\nu\partial_yu^\nu+\partial_xp^\nu-\nu\triangle u^\nu&=0,\\
 \partial_tv^{\nu}+u^\nu \partial_xv^\nu+v^\nu\partial_yv^\nu+\partial_yp^\nu-\nu\triangle v^\nu&=0,\\
 \partial_xu^\nu+\partial_yv^\nu&=0,\\
 (u^\nu,\;v^\nu)|_{y=0}&=(0,0),
\end{split}
\end{equation}
 where $(x,y)\in\mathbb{R}\times\mathbb{R}^+,$ $u^\nu$ is the tangential components of velocity to the
 boundary $(x,0)$, and $v^\nu$ is the normal components. A natural question is that does the solution
 $(u^\nu, v^\nu)$ convergence to the solution of Euler equation:
\begin{equation}
\begin{split}
 \partial_tu^E+u^E \partial_xu^E+v^E\partial_yu^E+\partial_xp^E&=0,\\
 \partial_tv^E+u^E \partial_xv^E+v^E\partial_yv^E+\partial_yp^E&=0,\\
 \partial_xu^E+\partial_yv^E&=0,\\
 v^\nu|_{y=0}&=0,
\end{split}
\end{equation}
 when $\nu\rightarrow0.$ As there is the no-slip condition: $u^\nu |_{y=0}=0$ in Navier-Stokes equations, the transition from zero
 velocity at the boundary to the full magnitude at some distance from it take place in a very thin layer.
 Then the flow can be divided into two regions: the boundary layer where viscous friction plays an essential
 part, and the remaining region outside this layer where friction may be neglected(the outer Euler flow). It was seen from
 several exact solution of Navier-Stokes equations that the boundary-layer thickness is proportional to $\sqrt{\nu}$, therefore
 we may write $(u^\nu,\;v^\nu)$ formally as:
\begin{eqnarray*}
u^\nu(t,x,y)&=&u^E(t,x,y)+u^B(t,x,\frac{y}{\sqrt{\nu}}),\\
v^\nu(t,x,y)&=&v^E(t,x,y)+\sqrt{\nu}v^B(t,x,\frac{y}{\sqrt{\nu}}),
\end{eqnarray*}
corresponding to the $p^\nu(t,x,y)$:
$$p^\nu(t,x,y)=p^E(t,x,y)+p^B(t,x,\frac{y}{\sqrt{\nu}}).$$
Denote $Y=\frac{y}{\sqrt{\nu}},$ let us define:
\begin{eqnarray*}
u(t,x,Y)&:=&u^E(t,x,0)+u^B(t,x,Y),\\
v(t,x,Y)&:=&Y\partial_y v^E(t,x,0)+v^B(t,x,Y),
\end{eqnarray*}
we formally obtain, from the Navier-Stokes equation by making $\nu$
tend to zero, the following system:

\begin{equation}
\begin{split}
\partial_tu+u\partial_xu+v\partial_Yu-\partial^2_Yu+\partial_x P=&0,\\
\partial_xu+\partial_Yv=&0,\\
(u,v)|_{Y=0}=&0,\\
\lim_{Y\rightarrow\infty}u=&u^E(t,x,0),
\end{split}
\end{equation}
where the pressure $p$ does not depend on $Y$, and satisfies the Bernoulli equation:
$$\partial_tu^E(t,x,0)+u^E(t,x,0)\partial_xu^E(t,x,0)+\partial_x p=0.$$

These are the Prandtl equations, proposed by Ludwig Prandtl \cite{p} in 1904. Although
the Prandtl equations have been simplified to a great extent, as compared with the Navier-Stokes
equations, they are still so difficult from the mathematical point of view that not very many
general statements about them can be made.

For the steady version, the von Mises transformation reduces these equations to a degenerated
parabolic equation under the condition $u>0.$
Let us consider the equation
in the domain $D=\{0<x<X,\;0<Y<\infty\}$ with the condition $$u(0,Y)=u_1(Y),$$
and introduce new independent variables by
$$\xi=x,\;\;\;\psi=\psi(x,Y),$$
where
$$u=\partial_Y\psi,\;\;\;v=\partial_x\psi,\;\;\;\psi(x,0)=0,$$
and a new function $\omega(\xi,\psi)=u^2(x,Y).$ Then the domain $D$ turns into $\{0<\xi<X,\;0<\psi<\infty\},$
and the Prandtl system reduces to the equation:
$$\sqrt{\omega}\partial^2_{\psi}\omega-\partial_\xi\omega-2p_{\xi}=0,$$
with the condition:
$$\omega(\xi,0)=0,\;\;\;\omega(0,\psi)=\omega_1(\psi),\;\;\;\lim_{\psi\rightarrow\infty}\omega(\xi,\psi)=U^2(x)$$
where
$$\omega_1(\int^Y_0u_1(\eta)d\eta)=u^2_1(Y).$$
Owing to this transformation, the maximum and
comparison principles apply for the equation, and we get the existence and uniqueness of (1.3):
\begin{Proposition}
(\cite{os})
Assume that $u_1(y)>0$ for $y>0;$ $u_1(0)=0,\,u'_1(0)>0,\,u_1(Y)\rightarrow U(0)\neq0$ as $Y\rightarrow\infty;$
$p_x$ is continuously differentiable on $[0,X];$ $u_1(Y),\,u'_1(Y),$ $u''_1(Y)$ are H\"{o}lder continuous and
 bounded for $0\leq Y<\infty.$ Moreover, assume that for small $y$ the compatibility condition is satisfied at
 the point $(0,0)$:
 $$u''_1(Y)-p_x(0)=O(Y^2).$$
 Then, for some $X>0$ there exists a unique solution $u(x,Y),\,v(x,Y)$ of the Prandtl equation (1.3) in $D$,
which have the following properties:

(1) $u(x,Y)$ is bounded and continuous in $\overline{D},$ $u>0$ for $Y>0;$

(2) $\partial_Y u>m>0$ for $0<Y\leq Y_0,$ where $m$ and $Y_0$ are constants;

(3) $\partial_Y u$ and $\partial^2_Y u$ are bounded and continuous in $D$;

(4) $\partial_xu,\,v$ and $\partial_Y v$ are bounded and continuous in any finite portion of $\overline{D}.$

Moreover,if $|u'_1(Y)|\leq m_1 e^{-m_Y},$ where $m_1$ and $m_2$ are positive constants, then
$\partial_xu$ and $\partial_Y v$ are bounded in $D$. If $p_x\leq 0,$ then such a solution exists in $D$
for any $X>0.$
\end{Proposition}

 We refer to \cite{os} for details. If
$p_x\leq 0,$ then $X=\infty,$ that is the separation of boundary layer does not appear under
this condition. When $p_x>0,$ the result is only valid local in $x,$ and this leads eventually
to boundary layer separation,see\cite{ca} \cite{e}.

The unsteady case is much complicated, and we have to imposed more conditions. Under the condition
 $\partial_yu>0,$ Crocco transformation, introduced by Crocco (1941), see \cite{c},
 reduces the boundary layer system to a single quasilinear equation of the degenerate parabolic
 type and satisfies the maximum principle. In contrast to the von Mises transformation, it reduces
 the Prandtl system to an equation in a finite domain and the boundary condition become nonlinear.
 basing on this transformation, Oleinik and Samokhin \cite{os} proved local in time
  well-posedness, and Xin and Zhang
 \cite{xz} proved the global in time well-posedness. Other positive mathematical results concern
  the case of analytic data. In \cite{sc}
  \cite{sc1} Sammartino,M.,and Caflisch,R.E.
  proved the short time existence and
uniqueness when the data are analytic; this result was improved
in \cite{ca} \cite{l} where analyticity is required in the x-axis
direction, while, using the regularizing effect of the viscosity.

In a recent result of G\'{e}rard-Varet and Dormy \cite{gd},  they established
ill-posedness for the Cauchy problem of the linearized Prandtl equation around non-monotic
special solution which is independent of $x$ and satisfies the heat equation. In \cite{gd2} and
\cite{gt}, some nonlinear ill-posedness were established with this counterexample.  As the system is global
well posed under the condition $\partial_yu>0.$ Then it is
natural to consider the problem that does this linear ill-posedness happen whenever the  non-degenerate
critical points appear.
As the well-posedness of stationary Prandtl equation has been established in \cite{os} and
the stationary solution could be considered as special time-periodic solution, and the
 time-periodic boundary layer system
\begin{equation}
\begin{split}
\partial_t u+u \partial_x u +v\partial_y u +p_x&=\partial^2_y u,\\
\partial_x u+\partial_y v&=0
\end{split}
\end{equation}
where $(t,x,y)\in\mathbb{T}\times[0,X)\times\mathbb{R}^+,$  with the condition
\begin{equation}
\begin{split}
u(t,x,0)=v(t,x,0)=0,\;\;u(t,0,y)=u_1(t,y),\;\;
\lim_{y\rightarrow\infty}u(t,x,y)=U(t,x),
\end{split}
\end{equation}
is well-posed under the condition $\partial_yu>0,$ see \cite{os},
then it is interesting to consider the linearized equation around stationary solution
$(u_0(x,y),v_0(x,y))$ where $u_0$ has non-degenerate
critical points, that is
\begin{equation}
\begin{split}
\partial_t u+u_0 \partial_x u +v_0\partial_y u +u\partial_xu_{0}+v\partial_y u_0-\partial^2_y u&=0,\\
\partial_x u+\partial_y v&=0
\end{split}
\end{equation}
where $(t,x,y)\in\mathbb{T}\times[0,X)\times\mathbb{R}^+,$  with the condition
\begin{equation}
\begin{split}
u(t,x,0)=v(t,x,0)=0,\;\;u(t,0,y)=u_1(t,y),\;\;
\lim_{y\rightarrow\infty}u(t,x,y)=0.
\end{split}
\end{equation}

Let us introduce the following function spaces:
$$W^{s,\infty}_{\alpha}:=\{f=f(y),\;e^{\alpha y}f\in W^{s,\infty}(\mathbb{R}^+)\},$$
with the norm:$\;\;\|f\|_{W^{s,\infty}_{\alpha}}=\|e^{\alpha y}f\|_{W^{s,\infty}},$
$$H^m_{\beta}:=H^m(\mathbb{T}_t,W^{0,\infty}_{\beta}(\mathbb{R}^+_y)),$$
and$$\overline{H}^m_{\beta}:=H^m_{\beta}\cap C^1(\mathbb{T}_t,W^{2,\infty}_{\beta}(\mathbb{R}^+_y)).$$
Our main result then reads:
\begin{theorem}
Let $u_0-U\in C^0([0,X_0);W^{4,\infty}_{\alpha}(\mathbb{R}^+))\cap C^1([0,X_0);W^{2,\infty}_{\alpha}(\mathbb{R}^+)),$
$u_0\mid_{x=0}$ has a non-degenerate critical point. If there exists $X>0,$ such that for every
$u_1(t,y)\in\overline{H}^m_{\beta},$ with $\beta<\alpha,$ equations (1.6) (1.7) have a unique solution,
let us denote $u(x,\cdot):=\mathfrak{X}(x,\xi)u_1,$ where $u(x,\cdot)$ is the solution of (1.6) and (1.7) with
 $u\mid_{x=\xi}=u_1,$ then there exist a $\delta>0,$ such that for every $\epsilon>0,$
 $$\sup_{0\leq\xi\leq x\leq \epsilon}\parallel e^{-\delta(x-\xi)\sqrt{|\partial_t|}}
 \mathfrak{X}(x,\xi)\parallel_{\mathcal{L}(H^m_\beta,H^{m-\sigma}_\beta)}=\infty,\;\;\forall m\geq 0,\;\sigma\in[0,\frac{1}{2}).$$
 If $\lim_{y\rightarrow\infty}u_0(x,y)=C,$ where $C\geq0$ is a constant, the result is valid for $\alpha=\beta.$
\end{theorem}

By modifying the construction of approximate solution performed in \cite{gd},  we construct an unstable quasimode which
is based on an asymptotic analysis of (1.6) in the high
time frequency limit and get the ill-posedness. We also need the following lemma:

\begin{lemma}
$C>0$ is a real constant, there exists $\tau\in\mathbb{C},$ with $im\tau<0$,
and a solution $W=W(z)$ of $$(\tau-z^2)^2\frac{d}{dz}W+iC\frac{d^3}{dz^3}((\tau-z^2)W)=0,$$
such that $$\lim_{z\rightarrow -\infty}W=0,\;\;\lim_{z\rightarrow +\infty}=1.$$
\end{lemma}

The only difference between this lemma and the the "spectral condition" in \cite{gd}
 is that there is a constance $C$ in this lemma. And this lemma could be proved by considering
 the eigenvalue problem:
 $$\frac{1}{z^2+1}u''+\frac{6z}{(z^2+1)^2}u'+\frac{6}{(z^+1)^2}u=\frac{\alpha}{C}u,$$
 for $z\in R,$ then extend $z$ from $R$ to $\mathbb{C}$, by the theory of ordinary differential
 equation and the complex change of variable, we get the existence of $W$. We refer to \cite{gd}
for details of the proof and list some property of $W$:
$$|W(z)-1|\leq C'e^{-c|z|^2},$$
when $z>0$,
$$|W(z)|\leq C'e^{-c|z|^2},$$
when $z<0,$ and
$$W^{(k)}(z)=O(e^{-c|z|^2}),\;\;\;\;\;z\rightarrow\infty,$$
where the constants $C',\;\;c>0.$

\begin{rmk}
Actually, the time-periodic ill-posedness of the equation correspondent to the large-time instability
of the equation. As in Z.Xin and L.Zhang \cite{xz}, it is well posed in global time with the condition $\partial_yu>0.$
Therefore there is some essential difference when the points satisfying $\partial_yu=0,\;\partial_yu\neq 0$ and depending
on $x$ appear.
\end{rmk}

\section{The Proof of The Result}
Firstly, let us construct the approximate solution by modifying the construction of D.G\'{e}rard-Varet and E.Dormy.
As equation (1.6) has constant coefficients in $t$, a Fourier analysis can be performed, then we could look for
solution in form $$u(t,x,y)=e^{-ikt}\hat{u}^k,\;\;v(t,x,y)=e^{-ikt}\hat{v}^k.$$
We can also separate corresponding frequency and amplitude in $x$ from $\hat{u}^k,$ which lead the ill-posedness.
That is we can look for the approximate solution in the form:
$$u_k(t,x,y)=e^{-it-i\omega(k)x}u^k(x,y).$$
We denote $\varepsilon:=\frac{1}{k}$ and concern the case $\varepsilon\ll 1.$

Denote $a$ the non-degenerate critical point of $u_0\mid_{x=0}.$ Assume that $\partial^2_y u_0(0,a)<0,$ then there
exists $0<X_1<X_0,\;a(x)\in C^1([0,X_1)),$ such that $$\partial_yu_0(x,a(x))=0,\;\partial^2_y(x,a(x))<0.$$ In fact,
$a(x)$ satisfies the differential equation:
$$\partial_{xy}u_0(x,a(x))+\partial^2_yu_0(x,a(x))a'(x)=0,$$
with the condition $a(0)=a.$

Let
$$u_{\varepsilon}=-ie^{-i\frac{t}{\varepsilon}-i\frac{1}{\varepsilon}\int^x_0\omega(\varepsilon,\xi)d\xi}
\frac{\partial_y(v^{reg}_{\varepsilon}+v^{sl}_{\varepsilon})}{\omega(\varepsilon,x)},\eqno(2.1)$$
$$v_\varepsilon=\frac{1}{\varepsilon}e^{-i\frac{t}{\varepsilon}-i\frac{1}{\varepsilon}\int^x_0\omega(\varepsilon,\xi)d\xi}
(v^{reg}_{\varepsilon}+v^{sl}_{\varepsilon})+ie^{-i\frac{t}
{\varepsilon}-i\frac{1}{\varepsilon}\int^x_0\omega(\varepsilon,\xi)d\xi}
(\frac{v^{reg}_{\varepsilon}+v^{sl}_{\varepsilon}}{\omega(\varepsilon,
x)})_x,\eqno(2.2)$$
where$$\omega(\varepsilon,x)=\frac{1}{-u_0(x,a(x))+\sqrt{\frac{\varepsilon}
{2}}\mid\partial^2_yu_0(x,a(x))\mid^{\frac{1}{2}}\tau},$$
the "regular" velocity:
$$v^{reg}_{\varepsilon}(x,y):=H(y-a(x))\big[u_0(x,y)-u_0(x,a(x))+
\sqrt{\frac{\varepsilon}{2}}\mid\partial^2_yu_0(x,a(x))\mid^{\frac{1}{2}}\tau\big],\eqno(2.3)$$
the shear layer velocity:
$$v^{sl}_{\varepsilon}(x,y):=\varphi(y-a(x))\sqrt{\frac{\varepsilon}{2}}
\mid\partial^2_yu_0(x,a(x))\mid^{\frac{1}{2}}
V[(\frac{\mid\partial^2_yu_0(x,a(x))\mid}{2\varepsilon})^{\frac{1}{4}}(y-a(x))],\eqno(2.4)$$
$\varphi$ is a smooth truncation function near 0,
$$V(z)=(\tau-z^2)(W(z)-H(z)),$$
$W$ satisfies the equation:
$$(\tau-z^2)^2\frac{d}{dz}W+iu_0(0,a)\frac{d^3}{dz^3}((\tau-z^2)W)=0,\eqno(2.5)$$
with the condition $$\lim_{z\rightarrow -\infty}W=0,\;\;\lim_{z\rightarrow +\infty}=1,$$
and $H$ is the Heaviside function. In the expression of $(u_\varepsilon,v_\varepsilon)$, the
"regular" velocity $v^{reg}_\varepsilon$ is the main part of the approximate solution which satisfies
the condition $$u_\varepsilon|_{y=0}=0,$$ and the shear layer velocity $v^{sl}_\varepsilon$ is the
modifying part as both $v^{reg}_\varepsilon$ and its second derivative have jumps at $y=a(x).$

It is easy to check that $u_\varepsilon$ is analytic in t, and $W^{2,\infty}_\beta$ in $x,y,$
and $$C_1e^{\frac{\delta_0t}{\sqrt{\varepsilon}}}\leq\|\frac{u_\varepsilon (t,x,y)}
{e^{-i\frac{t}{\varepsilon}}}\|_{W^{2,\infty}_\beta}\leq
C_2\frac{1}{\varepsilon^{\frac{1}{4}}}e^{\frac{\delta_0t}{\sqrt{\varepsilon}}},\eqno(2.6)$$
where $C_1,\;C_2$ and $\delta_0$ are independent of $\varepsilon.$

Inserting $u_\varepsilon,\,v_\varepsilon$ into (1.6), we have:
\begin{equation*}
\begin{split}
\partial_t u_\varepsilon+u_0 \partial_x u_\varepsilon +v_0\partial_y u_\varepsilon +u_\varepsilon
\partial_xu_{0}+v_\varepsilon\partial_y u_0-\partial^2_y u_\varepsilon&=I_\varepsilon,\\
\partial_x u_\varepsilon+\partial_y v_\varepsilon&=0,
\end{split}
\end{equation*}
moreover, $u_\varepsilon,\,v_\varepsilon$ satisfy the condition:
$$(u_\varepsilon,\,v_\varepsilon)\mid_{y=0}=(0,0),\;\lim_{y\rightarrow\infty}u_\varepsilon=0,$$
where
\begin{equation*}
\begin{split}
I_\varepsilon=&e^{-i\frac{t}{\varepsilon}-i\frac{1}{\varepsilon}\int^x_0\omega(\varepsilon,\xi)d\xi}
\big[-\frac{1}{\varepsilon}\frac{\partial_y(v^{reg}_{\varepsilon}+v^{sl}_{\varepsilon})}{\omega(\varepsilon,x)}
-\frac{u_0}{\varepsilon}\partial_y(v^{reg}_{\varepsilon}+v^{sl}_{\varepsilon})-iu_0\partial_x
\big(\frac{\partial_y(v^{reg}_{\varepsilon}+v^{sl}_{\varepsilon}
)}{\omega(\varepsilon,x)}\big)\\
&-iv_0\frac{\partial^2_y(v^{reg}_{\varepsilon}+v^{sl}_{\varepsilon})}
{\omega(\varepsilon,x)}-i\partial_x u_0\frac{\partial_y(v^{reg}_{\varepsilon}+v^{sl}_{\varepsilon})}
{\omega(\varepsilon,x)}+\frac{1}{\varepsilon}\partial_y u_0(v^{reg}_{\varepsilon}+v^{sl}_{\varepsilon})\\&
+i\partial_yu_0\partial_x\big(\frac{v^{reg}_{\varepsilon}+v^{sl}_{\varepsilon}}{\omega(\varepsilon,x)}\big)
+i\frac{\partial^3_y(v^{reg}_{\varepsilon}+v^{sl}_{\varepsilon})}{\omega(\varepsilon,x)}\big]
%e^{-i\frac{t}{\varepsilon}-i\frac{1}{\varepsilon}\int^x_0\omega(\varepsilon,\xi)d\xi}
\end{split}
\end{equation*}

For $y\neq a(x),\;I_\varepsilon$ could be write as
\begin{equation*}
\begin{split}
I_\varepsilon=e^{-i\frac{t}{\varepsilon}-i\frac{1}{\varepsilon}\int^x_0\omega(\varepsilon,\xi)d\xi}(I_1+I_2+I_3),
\end{split}
\end{equation*}
where
\begin{equation*}
\begin{split}
I_1=&-iu_0\partial_x
\big(\frac{\partial_y(v^{reg}_{\varepsilon}+v^{sl}_{\varepsilon}
)}{\omega(\varepsilon,x)}\big)-iv_0\frac{\partial^2_y(v^{reg}_{\varepsilon}+v^{sl}_{\varepsilon})}
{\omega(\varepsilon,x)}-i\partial_x u_0\frac{\partial_y(v^{reg}_{\varepsilon}+v^{sl}_{\varepsilon})}
{\omega(\varepsilon,x)}\\&+i\partial_yu_0\partial_x\big(\frac{v^{reg}_{\varepsilon}+v^{sl}_{\varepsilon}}
{\omega(\varepsilon,x)}\big),
\end{split}
\end{equation*}
\begin{equation*}
\begin{split}
I_2&=-\frac{1}{\varepsilon}\big(\frac{1}{\omega(\varepsilon,x)}+u_0(x,y)\big)\partial_yv^{reg}_\varepsilon
+\frac{1}{\varepsilon}\partial_yu_0v^{reg}_\varepsilon+i\frac{\partial^3_yv^{reg}_\varepsilon}{\omega(\varepsilon,x)}\\
&=i\frac{\partial^3_yv^{reg}_\varepsilon}{\omega(\varepsilon,x)},
\end{split}
\end{equation*}
and
$$I_3=-\frac{1}{\varepsilon}\big(\frac{1}{\omega(\varepsilon,x)}+u_0(x,y)\big)\partial_yv^{sl}_\varepsilon
+\frac{1}{\varepsilon}\partial_yu_0v^{sl}_\varepsilon+i\frac{\partial^3_yv^{sl}_\varepsilon}{\omega(\varepsilon,x)}.$$

As $u_0$ satisfies the equation
$$u_0\partial_xu_0-\int^y_0\partial_xu_0dy'\partial_yu_0+p_x=\partial^2_yu_0,$$
by solving the ordinary equation with respect to $y,$ $v_0(\cdot,y)$ could be write
as:
$$v_0(\cdot,y)=-u_0\int^y_0(\partial^2_0u_0-p_x)dy',$$
and
$$\partial_xu_0(\cdot,y)=\partial_yu_0\int^y_0(\partial^2_yu_0-p_x)dy'+u_0(\partial^2_yu_0-p_x),$$
then it is easy to check that
$$\partial_xu_0(\cdot,y)=O(1),\;\;\;v_0(\cdot,y)=O(y),\;\;\;as\;\;y\rightarrow\infty,$$
if $p_x\equiv0\,(U(x)\equiv C),$ we have
$$\partial_xu_0(\cdot,y)=O(e^{-\alpha y}),\;\;\;v_0(\cdot,y)=O(1),\;\;\;as\;\;y\rightarrow\infty.$$

Then by the property of $(u_0,\;v_0),$ we have
$$\|I_1(x,\cdot)\|_{W^{0,\infty}_\beta}\leq C_3,$$
where $C_3$ is a constant independent of $\varepsilon.$

As
\begin{eqnarray*}
I_3&=&-\frac{1}{\varepsilon}\big(\frac{1}{\omega(\varepsilon,x)}+u_0(x,y)\big)\partial_yv^{sl}_\varepsilon
+\frac{1}{\varepsilon}\partial_yu_0v^{sl}_\varepsilon+i\frac{\partial^3_yv^{sl}_\varepsilon}{\omega(\varepsilon,x)}\\
&=&-\frac{1}{\varepsilon}\big[u_0(x,y)-u_0(x,a(x))-\partial_yu_0(x,a(x))(y-a(x))\\
& &\;\;\;\;\;\;\;\;-\frac{1}{2}\partial^2_yu_0(x,a(x))(y-a(x))^2\big]\partial_yv^{sl}_\varepsilon\\
& &+\frac{1}{\varepsilon}\big[\partial_yu_0(x,y)-\partial_yu_0(x,a(x))-\partial^2_yu_0(x,a(x))
(y-a(x))\big]v^{sl}_\varepsilon\\
& &+\big\{-\frac{1}{\varepsilon}\big[\frac{1}{2}\partial^2_yu_0(x,a(x))(y-a(x))^2+\sqrt{\frac{\varepsilon}{2}}
|\partial^2_yu_0(x,a(x))|^{\frac{1}{2}}\tau\big]\partial v^{sl}_\varepsilon\\& &
\;\;\;\;\;\;+\frac{1}{\varepsilon}\partial^2_yu_0(x,a(x))(y-a(x))v^{sl}_\varepsilon-iu_0(0,a)
\partial^3_yv^{sl}_\varepsilon\big\}\\& &
+i\big[u_0(0,a)-u_0(x,a(x))+\sqrt{\frac{\varepsilon}{2}}|\partial^2_yu_0(x,a(x))|\big]\partial^3_yv^{sl}_\varepsilon\\
&:=&I_{31}+I_{32}+I_{33}+I_{34}.
\end{eqnarray*}
Set $z=(\frac{\mid\partial^2_yu_0(x,a(x))\mid}{2\varepsilon})^{\frac{1}{4}}(y-a(x)),$ it is easy to check
that:
$$\|I_{31}(x,\cdot)\|_{W^{0,\infty}_{\beta}}+\|I_{32}(x,\cdot)\|_{W^{0,\infty}_\beta}\leq C_4,$$
where $C_4$ is a constant independent of $\varepsilon$.
\begin{eqnarray*}
I_{33}&=&\frac{|\partial^2_yu_0(x,a(x))|^{\frac{1}{2}}}{(2\varepsilon)^{\frac{1}{2}}}
(z^2-\tau)\partial_yv^{sl}_\varepsilon-2(\frac{|\partial^2_yu_0(x,a(x))|}{2\varepsilon})^{\frac{3}{4}}
zv^{sl}_\varepsilon\\&=&(\frac{|\partial^2_yu_0(x,a(x))|}{2\varepsilon})^{\frac{3}{4}}
\sqrt{\frac{\varepsilon}{2}}|\partial^2_yu_0(x,a(x))|^{\frac{1}{2}}\varphi(y-a(x))
\big((z^2-\tau)V'-2zV\\& &-iu_0(0,a)\partial^3_yV'''\big)+\varphi'(y-a(x))(c_1\varepsilon^{k_1})V+c_2
\varepsilon^{k_2}V'')+\varphi''(y-a(x))c_3\varepsilon^{k_3}V'\\& &
+\varphi'''(y-a(x))c_4\varepsilon^{k_4}V,
\end{eqnarray*}
inserting the expression of $V$, we have:
$$(z^2-\tau)V'-2zV-iu_0(0,a)V'''=0,$$
when $z\neq0,$ and choosing $\varphi$, such that the derivative of $\varphi(y-a(x))$ is 0 near $y=a(x),$
as $V$ and its derivative decreases exponentially, then we have
$$\|I_{33}(x,\cdot)\|_{W^{0,\infty}_\beta}\leq C_5,$$
where $C_4$ is a constant independent of $\varepsilon.$

By the regularity of $u_0(x,y),$ we have
$$\|I_{34}(x,\cdot)\|_{W^{0,\infty}_\beta}\leq C_6,$$
when $0<x<\varepsilon^{\frac{1}{4}},$ where $C_6$ is a constant independent of $\varepsilon.$ Denote
$J_\varepsilon(x,y)=e^{i\frac{t}{\varepsilon}}I_\varepsilon(t,x,y),$ Then we have
$$\|J_\varepsilon(x,\cdot)\|_{W^{0,\infty}_\beta}\leq C_7e^{\frac{\delta_0x}{\sqrt{\varepsilon}}},\eqno(2.7)$$
when $0<x<\varepsilon^{\frac{1}{4}},$ where $C_7$ is a constant independent of $\varepsilon.$

Let us assume that (1.6) has a unique solution and $\forall\,\delta>0,
\;\exists\,\epsilon_0>0,\;m\geq0,$ and $0\leq\sigma<\frac{1}{2},$
$$\sup_{0\leq\xi\leq x\leq \epsilon_0}\parallel e^{-\delta(x-\xi)\sqrt{|\partial_t|}}
 \mathfrak{X}(x,\xi)\parallel_{\mathcal{L}(H^m,H^{m-\sigma})}\leq C_8,\eqno(2.8)$$

Let $\mathfrak{X}_\varepsilon(x,\xi)$ be the restriction of $\mathfrak{X}(x,\xi)$ to the
tangential Fourier mode $\frac{1}{\varepsilon}$, we have
$$\|\mathfrak{X}_\varepsilon(x,\xi)\|_{\mathcal{L}(W^{0,\infty}_\beta)}\leq
 C_8\varepsilon^{-\sigma}e^{\frac{\delta(x-\xi)}{\sqrt{\varepsilon}}}.\eqno(2.9)$$
 Let us write the first equation of (1.6) as
 \begin{eqnarray*}
 & &\partial_xu-\frac{\partial_xu_0}{u_0}\int^y_0 \partial_xudy'+\frac{\partial_tu}{u_0}
 +\frac{\partial_xu_0}{u_0}u+\frac{v_0}{u_0}\partial_yu-\partial^2_yu\\
 & :=&\partial_xu-\frac{\partial_xu_0}{u_0}\int^y_0 \partial_xudy'+Lu,
\end{eqnarray*}
and denote $L_\varepsilon$ the restriction of $L$ to the tangential Fourier mode $\frac{1}{\varepsilon}.$
Let $U(x,y)$ be the solution of $$\partial_xU-
\frac{\partial_xu_0}{u_0}\int^y_0\partial_xUdy'+L_\varepsilon U=0,$$
with the condition
$U(x,y)=e^{i\frac{t}{\varepsilon}}u_\varepsilon(t,x,y)$ when $x=0.$ Then
$$\|U(x,y)\|_{W^{0,\infty}_\beta}\leq C_9\varepsilon^{-\sigma}e^{\frac{\delta x}{\sqrt{\varepsilon}}}.\eqno(2.10)$$
Then $\widetilde{U}=U-e^{i\frac{t}{\varepsilon}}u_\varepsilon$ satisfies:
$$\partial_x\widetilde{U}-\frac{\partial_yu_0}{u_0}\int^y_0\partial_x\widetilde{U}dy'
+L_\varepsilon \widetilde{U}=J_\varepsilon,\eqno(2.11)$$
Assuming that the unique solution of (2.11) has the form
$\int^x_0\mathfrak{X}_\varepsilon(x,\xi)Q_\varepsilon(\xi,\cdot)d\xi,$ then we have
$$Q_\varepsilon(x,y)-\frac{\partial_yu_0}{u_0}\int^y_0Q_\varepsilon(x,y')dy'=J_\varepsilon(x,y),$$
by solving the differential equation we obtain
$$Q_\varepsilon(x,y)=J_\varepsilon(x,y)+\partial_yu_0\int^y_0J_\varepsilon(x,y')u_0(x,y')dy'.$$
It is easy to check that $\|Q_\varepsilon(x,\cdot)\|_{W^{0,\infty}_\beta}\leq C_{10},$ where $C_{10}$
 is a constant independent of $\varepsilon$. Then $\widetilde{U}=\int^x_0\mathfrak{X}\varepsilon(x,\xi)
 Q\varepsilon(x,\xi)d\xi,$ we have
 \begin{eqnarray*}
 \|\widetilde{U}(x,\cdot)\|_{W^{0,\infty}_\beta}&=&\sup_{y>0}|e^{\beta y}\int^x_0\mathfrak{X}_\varepsilon(x,\xi
 )Q_\varepsilon(\xi,y)d\xi|\\&\leq&\int^x_0\sup_{y>0}|e^{\beta y}\mathfrak{X}_\varepsilon(x,\xi
 )Q_\varepsilon(\xi,y)|d\xi\\&=&\int^x_0\|\mathfrak{X}_\varepsilon(x,\xi
 )Q_\varepsilon(\xi,y)\|_{W^{0,\infty}_\beta}d\xi\\&\leq&\int^x_0C_8\varepsilon^{-\sigma}e^
 {\frac{\delta(x-\xi)}{\sqrt{\varepsilon}}}\leq C_{11}\varepsilon^{\frac{1}{2}-\sigma}
 e^{\frac{\delta_0x}{\sqrt{\varepsilon}}},
 \end{eqnarray*}
 therefore
 $$
 \|U(x,\cdot)\|_{W^{0,\infty}_\beta}\geq \|e^{i\frac{x}{\varepsilon}}u_\varepsilon(x,\cdot)\|
 _{W^{0,\infty}_\beta}-\|\widetilde{U}(x,\cdot)\|_{W^{0,\infty}_\beta}\geq C_{12}e^
 {\frac{\delta_0x}{\sqrt{\varepsilon}}},\eqno(2.12)
 $$
 when $\delta<\delta_0.$ Then (2.12) contradicts (2.10), as soon as
 $\frac{\sigma}{\delta_0-\delta}|\ln\varepsilon|\sqrt{\varepsilon}\ll x<\varepsilon^{\frac{1}{4}}.$

\end{document}